# LES LOIS DE PROBABILITE
# POUR LES FONCTIONS STATISTIQUES

*(cas de collectifs à plusieurs dimensions)**

Par Y. Garti

(Institut des Mathématiques de l'Université d'Istanbul)

* Thèse de doctorat acceptée par la Faculté des Sciences de l'Université d'Istanbul.
(Rapporteurs: Prof. Dr. Richard von Mises, Prof. Dr. William Prager)

Remerciement

Nous adressons nos remerciements les plus sincères au Professeur Daniel Bennequin pour ses précieux conseils et son soutien.

## Les lois de probabilité pour les fonctions statistiques (cas de collectifs à plusieurs dimensions)

Dans un mémoire paru en 1936 dans les Annales de l'Institut Poincaré,[1] M. de Misès démontre que sous certaines conditions la distribution (loi de probabilité) des fonctions dites « statistiques » tend vers la Gaussienne, les collectifs étant supposés à une dimension. Voici en substance de quoi il s'agit.

Soit la suite infinie de collectifs à une dimension, de distributions respectives $V_1'(x), V_2'(x), \ldots$ et $x_1, x_2, \ldots, x_n$ le résultat d'une expérience effectuée sur les *n* premiers collectifs.

On sait depuis longtemps que la distribution de la somme de *n* variables aléatoires $x_1, x_2, \ldots, x_n : x_1 + x_2 + \ldots + x_n$ ou celle de leur moyenne arithmétique $\frac{1}{n}(x_1 + x_2 + \ldots + x_n)$ tend, pour *n* infini, vers la loi de Gauss sous certaines conditions très larges d'ailleurs.

M. de Misès montre que ce théorème s'étend à une classe beaucoup plus générale de fonctions des $x_1, x_2, \ldots, x_n$ : les fonctions statistiques. Donnons en une idée.

---

[1] Mises, R. von, 1936, Les lois de probabilité pour les fonctions statistiques. Annales de l'Institut Henri Poincaré, 6 (no. 3-4): pp. 185-212.

$S(x)$ étant la *répartition* des *n* variables $x_1, x_2, ..., x_n$ c'est-à-dire la fonction telle que, pour *x* donné, $nS(x)$ désigne le nombre de celles parmi les variables $x_1, x_2, ..., x_n$ dont la valeur ne surpasse pas *x*, toute fonction *f* des $x_1, x_2, ..., x_n$ qui s'exprime par $S(x)$ est une fonction statistique $f\{S(x)\}$.
Le résultat en question est:

Pour toute fonction statistique $f\{S(x)\}$ sous certaines conditions, la distribution tend vers la Gaussienne.

D'une façon plus précise, si l'on désigne par $S_n(x)$ la répartition des $x_1, x_2, ..., x_n$ et par $V_n(x)$ la moyenne

$$V_n(x) = \frac{1}{n}\left[V_1'(x) + V_2'(x) + ... + V_n'(x)\right]$$

des *n* premières distributions, dite distribution moyenne, sous certaines conditions, la probabilité $P_n(u)$ de l'inégalité

$$H_n\left[f\{S_n(x)\} - f\{V_n(x)\}\right] \leq u$$

où $H_n$ désigne une constante dépendant de *n*, converge uniformément pour toutes les valeurs de *u* vers la fonction de Gauss, si *n* tend vers l'infini:[2]

$$\lim_{n\to\infty} P_n(u) = \frac{1}{\sqrt{\pi}} \int^u e^{-t^2} dt = \phi(u).$$

Dans le présent exposé, je me propose de généraliser cette propriété dans le cas où les collectifs sont à plusieurs dimensions. Je ne transcrirai ici que les passages nécessitant une forme spéciale et principalement le cas des distributions continues, puisqu'en principe, pour les distributions arithmétiques, il n'y a pas de différence essentielle si la variable aléatoire, au lieu d'avoir une composante, en a plusieurs.

---

[2] *Le signe* $\int$ *désignera toujours l'integration de* $-\infty$ *à* $+\infty$ *, le signe* $\int_a$ *l'integration de* $a$ *à* $+\infty$ *etc.*

Citons quelques définitions spéciales au cas de plusieurs dimensions.

**Définitions**

Dans tout ce qui suit, $X$ désignera un point dans l'espace euclidien à $k$ dimensions $R_k$. En coordonnées rectangulaires, les composantes du point $X$ seront désignées par $x_1, x_2, ..., x_k, X = (x_1, x_2, ..., x_k)$. Les composantes de la variable aléatoire seront désignées par $\xi_1, \xi_2, ..., \xi_k$.

Les résultats possibles, c'est-à-dire l'espace caractéristique du collectif sera en général un domaine (sous ensemble) de $R_k$. Dans le cas où les résultats possibles forment un ensemble continu, la distribution sera continue, dans le cas contraire, c'est-à-dire que l'espace caractéristique n'est formé que de certains points de $R_k$, bien distincts, la distribution sera arithmétique ou discontinue.

Une fonction de distribution à $k$ dimensions est une fonction

$$V(X) = V(x_1, x_2, ..., x_k)$$

de $k$ variables, non négative et totalement monotone croissante, c'est-à-dire que:

**1°** si la distribution est continue et dérivable, on suppose la dérivée mixte

$\frac{\partial^k V}{\partial x_1 \partial x_2 ... \partial x_k}$ positive ou nulle:

$$\frac{\partial^k V}{\partial x_1 \partial x_2 ... \partial x_k} \geq 0.$$

**2°** si la distribution est arithmétique ou discontinue, *X* et *Y* étant deux points de $R_k$ tels que

$$x_i < y_i \qquad (i = 1, 2, ..., k)$$

l'expression correspondante à la dérivée précédente en différences finies en fonction de *X* et *Y*, est positive.

Ainsi dans le cas de deux variables $X = (x_1, x_2)$, $Y = (y_1, y_2)$ telles que $x_1 < y_1,\ x_2 < y_2$ on doit avoir

$$V(y_1, y_2) - V(y_1, x_2) - V(x_1, y_2) + V(x_1, x_2) \geq 0 \ .$$

$V(X)$ signifie la probabilité pour que les composantes $\xi_\nu$ de la variable ne dépassent pas $x_\nu$, c'est-à-dire la probabilité pour que

$$\xi_\nu \leq x_\nu \qquad (\nu = 1, 2, ..., k) \ .$$

La fonction $V(X)$ est continue à droite (du côté droit) c'est-à-dire que

$$V(x_1, x_2, ..., x_p, ..., x_k) = \lim_{\alpha_p \to 0} V(x_1, x_2, ..., x_p + \alpha_p, ..., x_k)$$

par valeurs positives $(p = 1, 2, ..., k)$.

De plus $V(X)$ prend respectivement les valeurs 0 et 1 à l'infini

$$V(-\infty, -\infty, ..., -\infty) = 0 \ , \ V(+\infty, +\infty, ..., +\infty) = 1 .$$

Soient $n$ points dans $R_k : X_1, X_2, ..., X_k$, par exemple les résultats de $n$ expériences effectuées sur $n$ collectifs quelconques. Désignons par $x_1^{(\mu)}, x_2^{(\mu)}, ..., x_k^{(\mu)}$ les composantes du point $X_\mu$ telles que

$$X_\mu = \left( x_1^{(\mu)}, x_2^{(\mu)}, ..., x_k^{(\mu)} \right) \ .$$

Nous définissons la « répartition » de ces $n$ résultats par une fonction de $k$ variables $S(X) = S(x_1, x_2, ..., x_k)$ telle que $nS(X)$ indique le nombre de résultats parmi les $X_\mu$, $(\mu = 1, 2, ..., n)$ dont les composantes sont inférieures ou égales à celles de $X$. Ces résultats sont ceux pour lesquels

$$\xi_i^{(\mu)} \le x_i \qquad (i = 1, 2, ..., k).$$

Ainsi dans le cas de deux dimensions $nS(X)$ désignera le nombre de points situés dans le rectangle construit sur $X$ comme diagonale et sur les côtés adjacents au point, c'est-à-dire le côté à droite et le côté supérieur.

Dans le cas de trois dimensions, $nS(x_1, x_2, x_3)$ désignera le nombre de points situés dans le parallélépipède construit sur $X$ comme diagonale et ceux qui sont aussi sur les trois faces adjacentes au point $X$, c'est-à-dire les trois faces qui se joignent au point $X$. On pourrait aussi généraliser cette interprétation géométrique. $S(X)$ est donc une fonction de *k* variables. Ses projections sur les axes donnent les répartitions à une dimension des coordonnées correspondantes de *X*.

Toute fonction qui s'exprime par $S(X)$ est une fonction statistique désignée par

$$f\{S(X)\} = f\{S(x_1, x_2, ..., x_k)\}.$$

Nous supposons toujours les distributions à $k$ variables.

Ainsi les moments

$$\int x_1^{\nu_1} x_2^{\nu_2} ... x_k^{\nu_k}\, dS(X)$$

sont des fonctions statistiques *linéaires,* et les écarts

$$M_{\nu_1 \nu_2 ... \nu_k} = \int (x_1 - a_1)^{\nu_1} (x_2 - a_2)^{\nu_2} \dots (x_k - a_k)^{\nu_k}\, dS(X)$$

avec

$$a_i = \int x_i\, dS(X) \qquad (i = 1, 2, ..., k)$$

sont des fonctions statistiques non linéaires.

De même toute fonction des moments, par exemple, le coefficient de corrélation dans le cas de deux dimensions

$$C=\frac{M_{11}}{\sqrt{M_{02}M_{20}}}$$

est une fonction statistique.

Nous montrerons que, sous certaines conditions, les lois de probabilité de ces fonctions tendent vers la Gaussienne.

## §. 1. Lemme préliminaire

**1.1**. Nous transcrirons ici un lemme tiré du mémoire de M. de Misès, ci-haut cité, en lui donnant une forme propre à notre cas.

Soient les collectifs $C_1, C_2, \ldots,$ à distributions quelconques. Désignons par $\alpha_n$ le nombre de dimensions de $C_n$, par $\xi_1, \xi_2, \ldots, \xi_{\alpha_n}$ les composantes de la variable aléatoire, par $V\left(\xi_1, \xi_2, \ldots, \xi_{\alpha_n}\right)$ la distribution dans $C_n$. En particulier $V_n\left(x_1, x_2, \ldots, x_{\alpha_n}\right)$ signifie la probabilité pour que dans le collectif $C_n$, $\xi_\nu$ ne dépasse pas $x_\nu : \xi_\nu \leq x_\nu$, $\left(\nu = 1, 2, \ldots, \alpha_n\right)$.

Soient deux fonctions quelconques des $x_\nu : A_n$ et $B_n$; les deux intégrales de Stieltjes suivantes, étendues à l'espace entier:

$$\int A_n \, dV_n = E\left\{A_n\right\}, \quad \int B_n \, dV_n = E\left\{B_n\right\}$$

sont les espérances mathématiques de $A_n$ et $B_n$ respectivement.

Les intégrales suivantes définissent les distributions de $A_n$ et $B_n$.

$$P_n(u) = \int_{A_n \leq u} dV_n, \qquad Q_n(u) = \int_{B_n \leq u} dV_n$$

La première de ces intégrales est étendue aux points de l'espace caractéristique pour lesquels la valeur de $A_n$ ne surpasse pas $u$; elle donne la probabilité de l'inégalité $A_n \leq u$.

Il en est de même pour la seconde intégrale.

Ces définitions faites, on démontre facilement le lemme suivant :

**1.2. Lemme** : *Soient* $A_n$ *et* $B_n$ *deux fonctions du caractère distinctif (de la variable aléatoire) dans le collectif* $C_n\left(A_n = A_n\left(X\right), B_n = B_n\left(X\right)\right)$ *dans une suite de collectifs* $C_1, C_2, \ldots$*; soient* $P_n\left(u\right)$ *et* $Q_n\left(u\right)$ *leurs distributions respectives, enfin*

$$\int \left|A_n - B_n\right| dV_n = E\left\{\left|A_n - B_n\right|\right\} = E_n$$

*l'espérance mathématique de leur différence en valeur absolue; les équations*

$$\lim_{n\to\infty} E_n = 0, \quad \lim_{n\to\infty} Q_n\left(u\right) = F\left(u\right)$$

*où* $F\left(u\right)$ *est une fonction à dérivée bornée:*

$$|F'(u)| < M$$

*entraînant l'équation*

$$\lim_{n\to\infty} P_n\left(u\right) = F\left(u\right)$$

## §. 2. Quelques formules concernant les épreuves répétées.

**2.1**. Les formules tirées aussi du mémoire de M. de Misès ont été adaptées ici au cas de collectifs à plusieurs dimensions.

Soit une suite infinie de collectifs $C'_1, C'_2, \ldots$ dont les caractères distinctifs ne sont que $l$ points de l'espace $R_k$ connus à l'avance et que nous caractériseront en attachant à chacun d'eux une valeur $f_1, f_2, \ldots, f_l$.

La probabilité pour que dans $C'_\nu$ la variable prenne la valeur $f_\lambda$ sera désignée par $p'_{\nu\lambda}$. Dès lors nous n'avons que des probabilités ponctuelles. M. de Misès démontre les propriétés suivantes:

Soit $C_n$ le collectif composé, obtenu en faisant une expérience sur les $n$ collectifs $C'_1, C'_2, ..., C'_n$.

Dans $C_n$ la probabilité d'un résultat tel que

$$f_{\alpha_1} f_{\alpha_2} \cdots f_{\alpha_n} \qquad (\alpha_1, \alpha_2, ..., \alpha_n = 1, 2, ..., l)$$

est

$$p'_{1\alpha_1} p'_{2\alpha_2} \cdots p'_{n\alpha_n}.$$

Envisageons dans $C_n$ des résultats différents quant à l'ordre des $f_\nu$ mais tels que la valeur $f_1$ s'y trouve $\rho_1 n$ fois, la valeur $f_2$, $\rho_2 n$ fois, … Appelons $P_n(\rho_1, \rho_2, ...\rho_l) = P_n(\rho)$ la probabilité de tels résultats (égale à la somme des valeurs $p'_{1\alpha_1} p'_{2\alpha_2} \cdots p'_{n\alpha_n}$ ).

On démontre que pour l'espérance mathématique, on a

$$\sum \rho_\lambda P_n(\rho) = E_n\{\rho_\lambda\} = p_{n\lambda}$$

avec

$$p_{n\lambda} = \frac{1}{n}\left( p'_{1\lambda} + p'_{2\lambda} + ... + p'_{n\lambda} \right) .$$

De même pour la dispersion, on a,

$$\sum \left[ \rho_\lambda - E_n\{\rho_\lambda\} \right]^2 P_n(\rho) = E_n\left\{ \left( \rho_\lambda - p_{n\lambda} \right)^2 \right\} \le \frac{1}{n} p_{n\lambda}(1 - p_{n\lambda}) \le \frac{1}{n} p_{n\lambda} .$$

Si l'on désigne par $\rho$ le vecteur de $l$ composantes $\rho_1, \rho_2, ...\rho_l$ et par $p_n$ le vecteur aux composantes $p_{n1}, p_{n2}, ..., p_{nl}$ on montre que

$$E_n\left\{ \left( \rho - p_n \right)^2 \right\} \le \frac{1}{n} .$$

**2. 2.** Jusqu'à présent nous avons supposé que les distributions dans les collectifs $C'_\nu$ fussent discontinues. Les résultats peuvent être étendus au cas de distributions continues.

Soit $V'_\nu(X) = V'_\nu(x_1, x_2, ...x_k)$ la probabilité pour que dans $C'_\nu$ la variable aléatoire ait ses composantes inférieures ou égales à celle de $X$ :

$$\xi_i \le x_i \qquad (i = 1, 2, ..., k).$$

Dans le collectif composé $C_n$ (par exemple résultat d'une expérience effectuée sur les $n$ premiers collectifs de la suite $C'_1, C'_2, ...,$) désignons par $nS(X)$ le nombre de résultats (de variables) dont les composantes ne surpassent pas celles de $X$. Soit $X = (x_1, x_2, ..., x_k)$ un vecteur bien déterminé, fixe. Dans $C'_\nu$, attribuons la valeur 1 à tous les résultats dont les composantes sont inférieures ou égales à celles de $X$ (points situés à l'intérieur ou sur le parallélépipède à $k$ dimensions, construit sur $X$ comme diagonale) et la valeur 0 à tous les autres résultats. Dans ce cas, $V'_\nu(X)$ et $1 - V'_\nu(X)$ désignent respectivement dans $C'_\nu$ les probabilités des résultats 1 et 0; $nS(X)$ est alors la somme de résultats qui se produisent dans une épreuve effectuée sur $C_n$. On peut alors appliquer la formule donnant l'espérance mathématique $E_n\{\rho_\lambda\}$ dans le cas des distributions discontinues, en posant $S(X)$ au lieu de $\rho_\lambda$ et $V'_\nu(X)$ au lieu de $p_{\nu\lambda}$:

$$E_n\{S(X)\} = \frac{1}{n}\left[V'_1(X) + V'_2(X) + ... + V'_n(X)\right];$$

et en posant

$$V_n(X) = \frac{1}{n}\left[V'_1(X) + V'_2(X) + ... + V'_n(X)\right]$$

on obtient pour les dispersions en appliquant la formule $y$ relative du cas arithmétique

$$E_n\left\{\left[S_n(X) - V_n(X)\right]^2\right\} \le \frac{1}{n} V_n(X)\left[1 - V_n(X)\right]$$

valable pour chaque point $X$.

Soit maintenant $\psi(X)$ une fonction non négative des composantes de $X : \psi(X) = \psi(x_1, x_2, ..., x_k)$ telle que l'intégrale

$$J = \int \psi(X)\left[S(X) - V_n(X)\right]^2 dX$$

existe. Il s'ensuit que:

$$E_n\{J\} =$$
$$= E_n\left\{ \int \psi(X)\left[S(X) - V_n(X)\right]^2 dX \right\} \le$$
$$\le \frac{1}{n}\int \psi(X) V_n(X)\left[1 - V_n(X)\right] dX$$

ou:

$$E_n\{J\} =$$
$$= E_n\left\{ \iint \ldots \int \psi\left(x_1, \ldots, x_p\right)\left[S\left(x_1, \ldots, x_k\right) - V_n(x_1, \ldots, x_k)\right]^2 dx_1 \ldots dx_k \right\} \le$$
$$\le \iint \ldots \int \psi\left(x_1, \ldots, x_p\right) V_n(x_1, \ldots, x_k)\left[1 - V_n\right] dx_1, \ldots, dx_k.$$

## §. 3. Cas de distributions arithmétiques.

**3.1**. Si les variables aléatoires ne peuvent prendre qu'une des $l$ valeurs fixées d'avance (valeurs attachées aux points résultats, c'est-à-dire si les résultats ne peuvent être qu'un des $l$ points de $R_k$) la répartition de $n$ résultats d'une épreuve composée est entièrement donnée par $l$ fréquences relatives $\rho_1, \rho_2, \ldots, \rho_l$ dont la somme est l'unité. Dans ce cas une fonction statistique des résultats est simplement une fonction de $l$ variables $\rho_1, \rho_2, \ldots, \rho_l$.

On peut partir d'une suite de collectifs $C_1', C_2', \ldots$ à distributions quelconques définies dans le même espace caractéristique ($R_k$ ou un de ses sous-ensembles). Par exemple on peut imaginer les tirages effectués dans des urnes dont chacune est remplie de billets portant chacun plusieurs chiffres. Nous subdivisons l'espace caractéristique en $l$ parties $L_1, L_2, ..., L_l$ (sous-ensembles de $R_k$ : par exemple par des parallélépipèdes généralisés).

Désignons par $p'_{\nu\lambda}$ la probabilité pour que le caractère distinctif de $C'_\nu$ tombe dans $L_\lambda$ $(\lambda = 1, 2, ..., l;\ \nu = 1, 2, ...)$. De cette façon les distributions primordiales sont réduites à des distributions discontinues, définies par les $p'_{\nu\lambda}$. On voit que le fait que les collectifs puissent être à plusieurs dimensions n'influe plus. La démonstration se déduit aisément comme dans le cas d'une dimension.

Nous citerons le théorème de limite concernant les distributions arithmétiques:

Soit $C_n$ la composition des collectifs $C'_1, C'_2, ..., C'_n$ (la suite des tirages effectués dans les *n* premières urnes). Soit $p_{n\lambda}$ la moyenne arithmétique des $p'_{\nu\lambda}$. Le résultat d'un tirage composé ou un élément de $C'_n$ comprend $\rho_1 n$ valeurs tombant dans $L_1$, $\rho_2 n$ dans $L_2$, $\rho_l n$ dans $L_l$.

Les variables $\rho_\lambda$ varient dans un domaine D défini par les relations:

$$\rho_1 \geq 0, \rho_2 \geq 0, ..., \rho_l \geq 0, \qquad \rho_1 + \rho_2 + ... + \rho_l = 1$$

(le simplex à $l$ dimensions).

Désignons par $\rho$ le vecteur aux composantes $\rho_\lambda$ et par $p_n$ celui aux composantes $p_{n\lambda}$. La dérivée d'une fonction $f(\rho) = f(\rho_1, ..., \rho_l)$ par rapport à $\rho_\lambda$ au point $p_n$ sera désignée par $f_\lambda$. On a le théorème suivant.

**3.2. Théorème**

*Soit* $f(\rho)$ *une fonction des fréquences relatives jouissant des propriétés suivantes :*

**1°** $f(\rho)$ *est borné dans le domaine D défini ci-dessus.*

**2°** $f(\rho)$ *admet des dérivées continues et bornées de premier et deuxième ordre.*

**3°** *Il existe deux indices* $\alpha$ *et* $\beta$ *différents l'un de l'autre et un nombre positif* $\eta$ *tels que à partir d'un certain* $n$

$$p'_{n\alpha} > \eta, \; p'_{n\beta} > \eta, \qquad \left| f_\alpha - f_\beta \right| > \eta.$$

*Sous ces conditions, la distribution des probabilités pour* $f(\rho)$ *tend vers la distribution de Gauss, si* $n$ *augmente de plus en plus. En d'autres termes, la probabilité* $P_n(u)$ *de l'inégalité*

$$H_n \left[ f(\rho) - f(p_n) \right] \leq u$$

*satisfait uniformément pour toutes les valeurs de* $u$ *à l'équation:*

$$\lim_{n \to \infty} P_n(u) = \phi(u) = \frac{1}{\sqrt{\pi}} \int^u e^{-t^2} dt$$

où $H_n$ est défini par

$$\frac{1}{2H_n^2} = \frac{1}{n^2} \sum_{\nu=1}^{n} \left[ \sum_{\lambda=1}^{l} f_\lambda^2 p'_{\nu\lambda} - \left( \sum_{\lambda=1}^{l} f_\lambda p'_{\nu\lambda} \right)^2 \right].$$

La démonstration utilise le fait que, sous certaines conditions, la distribution d'une somme de variables indépendantes tend vers la Gaussienne. Elle se base, en second lieu, sur le lemme.

**3. 3.** On montre aussi en remarque que la moyenne $b_n$ et la dispersion $s_n^2$ de la distribution limite de $f(\rho)$ sont respectivement

$$b_n \sim f(p_n) \qquad \text{et} \qquad s_n^2 \sim \frac{1}{n^2} \sum_{\nu=1}^{n} \left[ \sum_{\lambda=1}^{k} f_\lambda^2 p'_{\nu\lambda} - \left( \sum_{\lambda=1}^{k} f_\lambda p'_{\nu\lambda} \right)^2 \right].$$

## §. 4. Les fonctions statistiques.

**4.1**. Si on reprend la définition générale que nous avons donnée pour les distributions, la répartition de *n* résultats, $X_1, X_2, ..., X_n : S(X)$ est un cas particulier d'une distribution.

Nous avons à considérer certains ensembles $J$ de distributions comprenant des répartitions à valeur quelconque de *n* et des distributions d'autre nature. Soient $V(X)$ et $V_1(X)$ deux distributions de $J$.
Si l'on fait varier $t$ entre 0 et 1, les distributions

$$V_1(X) + t\left[V(X) - V_1(X)\right] \qquad 0 \le t \le 1$$

constituent le segment de droite de $V_1(X)$ à $V(X)$. Un ensemble qui comprend tous ces segments déterminés par deux de ses éléments s'appellera un ensemble « *convexe* ». Un tel ensemble convexe *J* pourra être défini par exemple, par toutes les distributions $V(X)$ pour lesquelles en un certain point $(X) = (X_1)$, la valeur $V(X_1)$ est comprise entre deux valeurs $a$ et $b$ positives et inferieures à l'unité. Ou bien un ensemble convexe est constitué par les distributions $V(X)$ pour lesquelles le produit

$$V(X)\left[1 - V(X)\right]$$

s'annule pour $|X|$ infini,

$$\left(|x_1| = \infty, |x_2| = \infty, ..., |x_k| = \infty\right)$$

comme une certaine puissance négative de $|X| = \sqrt{x_1^2 + x_2^2 + ... + x_k^2}$ .

**4.2**. Soit $J$ un ensemble de distributions suivant les explications données. Attachons à chaque fonction $V(X)$ de $J$, une certaine valeur *f*. *f* est une *fonction statistique* définie sur $J$ et on écrit $f\{V(X)\}$.

Les fonctions dont on s'occupe dans la statistique générale, moyennes, écarts quadratiques, moments, coefficients de corrélation, quotient de Lexis,…sont des fonctions statistiques. Si dans l'expression de $f\{V(X)\}$ on substitue à $V(X)$ une répartition $S(X)$, $f$ devient une fonction de $n$ points $X_1, X_2, ..., X_n$ fonction jouissant des propriétés suivantes:

**1)** Elle est symétrique, c'est-à-dire que la valeur de $f$ ne change pas si on remplace simultanément le résultat $X_\mu$ par $X_\lambda$ et $X_\lambda$ par $X_\mu$.

**2)** La valeur de $f$ ne change pas si l'on passe de $n$ à un multiple entier $2n, 3n, \ldots$ et que tout résultat $X_1, X_2, ..., X_n$ se retrouve 2 fois, 3 fois,… parmi les nouvelles variables.

Si l'on admet comme variable indépendante de $f$ que des fonctions $V(X)$ et $S(X)$ qui restent constantes partout, sauf en certains points connus au préalable $A_1, A_2, ..., A_m$, $f$ s'exprime comme une fonction de $m$ variables ordinaires, à savoir comme fonction des fréquences relatives $\rho_1, \rho_2, ..., \rho_m$ dans le cas de $S$ et comme fonction de probabilités ponctuelles $p_1, p_2, ..., p_m$ au cas de $V$. Il ne s'agit plus dès lors que de distributions arithmétiques et répartitions à points de saut fixes, problème déjà traité.

**4.3.** Comme exemple de fonction statistique *linéaire* on peut citer l'intégrale de Stieltjes $f\{V(X)\} = \int \psi(X)\, dV$, $\psi(X)$ étant une fonction continue de plusieurs variables. Dans ce genre on peut citer les moments

$$\int x_1^{\nu_1}, x_2^{\nu_2}, ..., x_k^{\nu_k}\, dV(X).$$

Si l'on prend pour $V(X)$ la répartition $S(X)$ des résultats $X_1, ..., X_n$, on a

$$f\{S(X)\} = \frac{1}{n}\sum_{\nu=1}^{n}\psi(X_\nu) = \frac{1}{n}\left[\psi(X_1)+...+\psi(X_n)\right].$$

Si $\psi$ est borné, $f$ est défini pour toute distribution $V(X)$; au cas contraire l'ensemble $J$ où $f$ existe, est limité par certaines conditions à remplir pour les $V(X)$ à l'infini. On supposera absolument convergentes les intégrales généralisées que l'on rencontrera dans les définitions des fonctions statistiques.

On peut citer comme exemple de fonctions *non linéaires* les écarts

$$M_{\nu_1\nu_{2.}...\nu_k} = \int (x_1-\alpha_1)^{\nu_1}(x_2-\alpha_2)^{\nu_2}\cdots(x_k-\alpha_k)^{\nu_k}\, dV(X)$$

avec

$$\alpha_i = \int x_i\, dV \qquad (i=1,2,...,k).$$

D'une façon générale toute fonction ressortant d'une combinaison d'intégrales de Stieltjes est une fonction statistique.

On peut citer d'autres cas aussi:

$$f\{V(X)\} = \iint \psi(X,Y)\, dV(X)\, dV(Y)$$

ou explicitement

$$f\{V(X)\} =$$
$$= \int...\int \psi(x_1,x_2,...,x_k;y_1,y_2,...,y_k)\, dV(x_1,x_2,...,x_k)\, dV(y_1,y_2,...,y_k)\ .$$

**4.4. Dérivées d'une fonction statistique.**

**4.4.1.** Soit $f\{V(X)\}$ définie sur un ensemble convexe $J$; soit $V_1(X)$ une distribution déterminée appartenant à $J$. On dit que $f\{V(X)\}$ est dérivable au point $V_1(X)$ si:

a) la fonction de $t$

$$F(t) = f\left\{V_1(X) + t(V - V_1)(X)\right\}$$

est dérivable par rapport à $t$ pour $t = 0$, quelle que soit la distribution $V(X)$ de $J$.

b) Cette dérivée qui dépendra de $V(X)$ et $V_1(X)$ s'exprime par une intégrale de Stieltjes.

$$\frac{d}{dt} f\left\{V_1(X) + t(V - V_1)(X)\right\}_{t=0} = \int f'\left\{V_1(X), Y\right\} d(V - V_1)(Y) \, .$$

($Y$ étant toujours une variable à $k$ composantes); $f'$ dépend de $V_1$ et de $Y$ *mais pas de* $V(X)$.

Ces conditions remplies, nous appellerons $f'\left\{V(X), Y\right\}$ la dérivée de la fonction statistique $f$ au point $V(X)$.

La dérivée de la fonction linéaire

$$f\left\{V(X)\right\} = \int \psi(X)\, dV(X)$$

ne dépend pas de $V_1$ et est tout simplement égale à $\psi(Y)$; en effet, on a:

$$\frac{d}{dt}\int \psi(X)\, d\left[V_1(X) + t(V - V_1)(X)\right]_{t=0} = \int \psi(Y)\, d(V - V_1)(Y)$$

et ainsi

$$f'\left\{V(X), Y\right\} = \psi(Y).$$

Si $f$ est une fonction de plusieurs intégrales de Stieltjes de la forme

$$A = \int \alpha(X)\, dV(X), \;\; B = \int \beta(X)\, dV(X), \;\; C = \int \gamma(X)\, dV(X), \;\; \ldots$$

c'est-à-dire que $f = F(A, B, C, \ldots)$,

la dérivée de la forme

$$f'\left\{V(X), Y\right\} = \frac{\partial F}{\partial A}\alpha(Y) + \frac{\partial F}{\partial B}\beta(Y) + \ldots$$

(les fonctions $\alpha, \beta, \ldots$peuvent dépendre d'un nombre quelconque de variables, mais toutefois ne dépassant pas *k*). Ceci se démontre facilement. Ainsi par exemple, pour l'écart

$$M_m^{(i)} = \int (x_i - a_i)^m \, dV(X) \text{ où } a_i = \int x_i \, dV$$

se rapportant à la composante $x_i$ de la variable aléatoire, on a

$$f'\{V(X), Y\} = (y_i - a_i)^m - mM_{m-1}^{(i)} y_i.$$

Considérons l'écart plus général

$$M_{\nu_1 \nu_2 \ldots \nu_k} = \int (x_1 - a_1)^{\nu_1} (x_2 - a_2)^{\nu_2} \ldots (x_k - a_k)^{\nu_k} \, dV$$

avec

$$\alpha_i = \int x_i \, dV \qquad (i = 1, 2, \ldots k)$$

et calculons sa dérivée.

Posons

$$A = \int (x_1 - a_1)^{\nu_1} \ldots (x_k - a_k)^{\nu_k} \, dV, \; a_1 = B = \int x_1 \, dV, \; a_2 = C = \int x_2 \, dV, \ldots$$

On a:

$$f'\{V(X), Y\} = \frac{\partial F}{\partial A} \alpha(Y) + \frac{\partial F}{\partial B} \beta(Y) + \ldots$$

Ici comme on le voit

$$\alpha = (x_1 - a_1)^{\nu_1} \ldots (x_k - a_k)^{\nu_k}, \; \beta = x_1, \; \gamma = x_2, \ldots$$

or:

$$\frac{\partial F}{\partial A} = 1 \,,$$

$$\frac{\partial F}{\partial B} = -\nu_1 \int (x_1 - a_1)^{\nu_1 - 1} (x_2 - a_2)^{\nu_2} \ldots (x_k - a_k)^{\nu_k} \, dV$$

$$= -\nu_1 M_{\nu_1 - 1, \nu_2, \ldots, \nu_k}$$

et ainsi de suite. D'où

$$f'\{V(X),Y\}=$$

$$=(y_1-a_1)^{\nu_1}(y_2-a_2)^{\nu_2}\ldots(y_k-a_k)^{\nu_k}$$
$$-\nu_1 M_{\nu_1-1,\nu_2,\ldots,\nu_k}y_1-\nu_2 M_{\nu_1,\nu_2-1,\nu_3,\ldots,\nu_k}y_2-\ldots-$$
$$-\nu_k M_{\nu_1,\ldots,\nu_k-1}y_k.$$

La fonction non linéaire

$$f=\iint\psi(X,Y)\,dV(X)\,dV(Y)$$

admet la dérivée

$$f'\{V(X),Y\}=\int[\psi(X,Y)+\psi(Y,Z)]dV(X)$$

comme on le verra ci-dessous.

***Remarque***

Une constante additive (expression ne dépendant pas de $Y$) ajoutée à $f'$ n'a aucune importance, étant donné que

$$\int d(V-V_1)$$

s'annule toujours.

**4.4.2.** Nous aurons aussi besoin dans ce qui suit de la « deuxième dérivée » d'une fonction statistique $f\{V(X)\}$. Cette dérivée sera désignée par $f''\{V(X),Y,Z\}$ où $Y$ et $Z$ sont des variables à plusieurs composantes et définie par la relation

$$\frac{d^2}{dt^2}f\{V_1(X)+t(V-V_1)(X)\}_{t=0}=$$
$$=\iint f''\{V_1(X),Y,Z\}\,d(V-V_1)(Y)\,d(V-V_1)(Z).$$

La seconde dérivée d'une fonction linéaire de la forme

$$f\{V(X)\} = \int \psi(X)\, dV(X)$$

est nulle.

Pour des fonctions de la forme $f = F(A, B, ...)$ dont on a parlé ci-dessus, on trouve

$$f''\{V(X), Y, Z\} = \frac{\partial^2 F}{\partial A^2} \alpha(Y)\ \alpha(Z) + 2\frac{\partial^2 F}{\partial A \partial B} \alpha(Y)\, \beta(Z) + ...$$

**4.4.3. Exemples**

a) L'écart $M_m^{(i)} = \int (x_i - a_i)^m\, dV(X)$ admet la seconde dérivée

$$f''\{V(X), Y, Z\} = -2m(y_i - a_i)^{m-1} z_i + m(m-1) M_{m-2}^{(i)} y_i z_i$$

que l'on peut déduire de l'écart général cité ci-dessous.

Pour $m = 2$ on a $f'' = -2y_i z_i$: on a supprimé le terme $4a_i z_i$ car une expression ne dépendant que d'un seul des points $Z$ ou $Y$ n'a pas d'importance sur $f''$ car $\int d(V - V_1)(X) = 0.$

b) Soit l'écart plus général

$$M_{\nu_1 \nu_2 ... \nu_k} = \int (x_1 - a_1)^{\nu_1} (x_2 - a_2)^{\nu_2} \cdots (x_k - a_k)^{\nu_k}\, dV$$

la seconde dérivée est égale à:

$$\begin{aligned} f''\{V(X), Y, Z\} &= \\ &= \nu_1(\nu_1 - 1) M_{\nu_1 - 2, \nu_2, ..., \nu_k} y_1 z_1 + \nu_2(\nu_2 - 1) M_{\nu_1, \nu_2 - 2, ..., \nu_k} y_2 z_2 + \\ &+ ... + \nu_k(\nu_k - 1) M_{\nu_1, \nu_2, ..., \nu_k - 2} y_k z_k - \\ &- 2\nu_1 (y_1 - a_1)^{\nu_1 - 1} (y_2 - a_2)^{\nu_2} ... (y_k - a_k)^{\nu_k} z_1 - \\ &- 2\nu_2 (y_1 - a_1)^{\nu_1} (y_2 - a_2)^{\nu_2 - 1} (y_3 - a_3)^{\nu_3} ... (y_k - a_k)^{\nu_k} z_2 - ... \\ &... - 2\nu_k (y_1 - a_1)^{\nu_1} ... (y_k - a_k)^{\nu_k - 1} z_k. \end{aligned}$$

Il est facile de le constater.

c) Considérons la fonction

$$f\left\{V(X)\right\}=\iint\psi(X,Y)\,dV(X)\,dV(Y)$$

On a:

$$F(t)=\iint\psi(X,Y)\,d\left[V_1(X)+t(V-V_1)(X)\right]\left[V_1(Y)+t(V-V_1)(Y)\right]=$$
$$=\iint\psi(X,Y)\Big[dV_1(X)\,dV_1(Y)+$$
$$+t\big(dV_1(X)\,d(V-V_1)(Y)+dV_1(Y)\,d(V-V_1)(X)\big)+$$
$$+t^2d(V-V_1)(X)\,d(V-V_1)(Y)\Big].$$

Prenons la première dérivée de $F$ par rapport à $t$ et faisons $-y$ $t=0$.

On obtient

$$F'(0)=\int\int\psi(X,Y)\,dV_1(X)\,d(V-V_1)(Y)+$$
$$+\int\int\psi(X,Y)\,dV_1(Y)\,d(V-V_1)(X)=$$
$$=\int\int\psi(X,Y)\,dV_1(X)\,d(V-V_1)(Y)+$$
$$+\int\int\psi(Y,X)\,dV_1(X)\,d(V-V_1)(Y).$$

On voit que:

$$f''\left\{V(X),Y\right\}=\int\left[\psi(X,Y)+\psi(Y,X)\right]dV(X).$$

Quant à la deuxième dérivée on obtient en calculant $F''(0)$

$$F''(0) =$$
$$= 2\iint \psi(X,Y)\, d(V-V_1)(X)\, d(V-V_1)(Y) =$$
$$= 2\iint \psi(Y,Z)\, d(V-V_1)(Y)\, d(V-V_1)(Z),$$

ce qui montre que

$$f''\{V(X),Y,Z\} = 2\,\psi(Y,Z).$$

**4.5. Formule de Taylor pour une fonction statistique.**

On considère la fonction

$$F(t) = f\left\{V_1(X) + t(V-V_1)(X)\right\};$$

on peut écrire le développement limité

$$F(1) - F(0) = F'(0) + \frac{1}{2}F''(\vartheta) \qquad 0 \le \vartheta \le 1$$

comme dans le cas d'une dimension. En utilisant la variable auxiliaire

$$t' = \frac{t-\vartheta}{1-\vartheta}$$

on déduit du développement ci-dessus en considérant la définition des dérivées de $f\left\{V(X)\right\}$, d'une façon identique au cas d'une dimension, le développement suivant:

$$f\left\{V(X)\right\} - f\left\{V_1(X)\right\} =$$
$$= \int f'\{V_1(X),Y\}\, d(V-V_1)(Y) +$$
$$+ \frac{1}{2}\iint f''\{V_2(X),Y,Z\}\, d(V-V_1)(Y)\, d(V-V_1)(Z).$$

İci $V_2$ désigne une distribution définie pour

$$V_2(X) = V_1(X) + \vartheta(V-V_1)(X)$$

pour $\vartheta$ compris entre 0 et 1.

Une telle formule subsiste pour toute fonction statistique deux fois dérivable sur le « segment » allant de $V_1(X)$ à $V(X)$.

***Remarque***

Il se peut que le premier et le second terme du développement soient identiquement nuls, comme cela a lieu pour la fonction

$$\omega^2\{V(X)\} = \int \lambda(X)\left[V(X)-V_1(X)\right]^2 dX$$

qu'il faut étudier spécialement.

## §. 5. Le théorème général

**5. 1.** On donne une suite infinie de collectifs $C_1', C_2'$ à $k$ dimensions. Soient respectivement

$$V_1'(X)=V_1'(x_1,x_2,...,x_k),\ V_2'(X)=V_2'(x_1,x_2,...,x_k),...$$

les distributions correspondantes supposées continues.

Une expérience effectuée sur les $n$ premiers collectifs fournit $n$ résultats $X_1, X_2,..., X_n$ (points dans $R_k$ ) dont la répartition sera désignée par $S_n(X)$.

Soit $f\{V(X)\}$ une fonction statistique ; l'expérience envisagée fait ressortir une certaine valeur $f\{S_n(X)\}$.

On sait d'après les règles élémentaires du calcul des probabilités que la valeur de $f\{S_n(X)\}$ est soumise à une certaine loi de probabilité, $P_n(u)$ entièrement définie par $V_1'(X), V_2'(X),..., V_n'(X)$.

Notamment $H_n$ et $K_n$ étant deux constantes, nous désignerons par $P_n(u)$ la probabilité pour que la valeur $H_n\left[f\{S_n(X)\}-K_n\right]$ déduite de l'expérience sur les $n$ premiers collectifs ne surpasse pas $u$ .

Nous nous proposons de démontrer que sous certaines conditions $P_n(u)$ tend vers la distribution de *Gauss* pour $n$ infini.

Il peut arriver que les répartitions qui se produisent au cours d'une expérience soient soumises à une certaine restriction, par exemple que les points de saut ne se trouvent que dans la partie de l'espace caractéristique définie par les coordonnées positives. En tout cas, il y aura un certain ensemble de répartitions admissibles pour toute valeur de $n$. D'autre part les distributions données donnent lieu à une autre suite de distributions, les moyennes arithmétiques définies par

$$V_n(X) = \frac{1}{n}\left[V_1'(X) + V_2'(X) + \ldots + V_n'(X)\right]$$

Dans ce qui suit, on supposera que $f\{V(X)\}$ soit définie dans un ensemble convexe *J* qui comprend toutes les répartitions admissibles et tous les $V_n(X)$ au moins à partir d'un certain $n$.

### 5.2. Théorème

*Soit* $f\{V(X)\}$ *une fonction statistique satisfaisant aux conditions suivantes:*

1) $f\{V(X)\}$ es*t deux fois dérivable dans un ensemble convexe J comprenant toutes les répartitions* $S_n(X)$ *qui peuvent se présenter au cours des expériences et toutes les distributions moyennes* $V_n(X)$ *au moins à partir d'un certain* $n$.

2) *La première dérivée* $f'\{V_n(X), Y\}$ *remplit des conditions suffisantes pour la validité du théorème classique sur l'application de la loi de Gauss.*

Par exemple si pour $\nu = 1, 2, 3, \ldots$

$$a_\nu = \int f'\{V_n(X),Y\}\, dV'_\nu(Y), \quad r_\nu^2 = \int\left[ f'\{V_n(X),Y\} - a_\nu \right]^2 dV'_\nu(Y),$$
$$C_\nu = \int \left| f'\{V_n(X),Y\} - a_\nu \right|^{2+\varepsilon} dV'_\nu(Y),$$

il suffit que $s_n^2 = r_1^2 + r_2^2 + \ldots + r_n^2$ divisé par $n^{\frac{1}{2+\varepsilon}}$ $(0 < \varepsilon < 1)$ tende vers l'infini et que les $C_\nu$ soient bornés.

3) *Pour la seconde dérivée* $f''\{V(X),Y,Z\}$ *il existe une fonction positive* $\psi(X)$ *telle que pour* $T_n(X) = V_n(X) - S_n(X)$ *quelle que soit la distribution* $V(X)$ *de J, l'inégalité*

$$\left| \iint f''\{V(X),Y,Z\}\, dT_n(Y)\, dT_n(Z) \right| \leq \int \psi(X)\, T_n^2 dX$$

*entraine la relation*

$$\lim_{n\to\infty} \frac{1}{s_n} \int \psi(X) V_n(X)\left[1 - V_n(X)\right] dX = 0.$$

*Sous ces conditions, la distribution des probabilités des* $f\{S_n(X)\}$ *tend vers la distribution de Gauss pour* $n$ *infini. En d'autres termes, si* $P_n(u)$ *désigne la probabilité de l'inégalité*

$$H_n\left[ f\{S_n(X)\} - f\{V_n(X)\} \right] \leq u$$

*on a uniformément pour toutes les valeurs de u,*

$$\lim_{n\to\infty} P_n(u) = \frac{1}{\sqrt{\pi}} \int^u e^{-t^2} dt = \phi(u),$$

*où* $H_n$ *est déterminé par*

$$\frac{1}{2H_n^2} = \frac{s_n^{\;2}}{n^2} =$$

$$= \frac{1}{n^2} \sum_{\nu=1}^{n} \left[ \int f'^2 \left\{ V_n(X), Y \right\} dV'_\nu(Y) - \left( \int f' \left\{ V_n(X), Y \right\} dV'_\nu(Y) \right)^2 \right].$$

**5.3.** ***Démonstration***

**5.3.1.** Elle se fait comme pour les distributions à une variable. Nous la transcrivons rapidement. Utilisons la formule de Taylor, établie précédemment

$$f\left\{V(X)\right\} - f\left\{V_1(X)\right\} =$$
$$= \int f'\left\{V_1(X), Y\right\} d\left(V - V_1\right)(Y) +$$
$$+ \frac{1}{2} \iint f''\left\{V_2(X), Y, Z\right\} d\left(V - V_1\right)(Y)\, d\left(V - V_1\right)(Z).$$

Son emploi est permis en vertu de la première hypothèse faite sur $f\left\{V(X)\right\}$. En prenant $S_n(X)$ pour $V(X)$ et $V_n(X)$ à la place de $V_1$ on obtient après avoir multiplié par $H_n$

$$H_n\left[ f\left\{S_n(X)\right\} - f\left\{V_n(X)\right\} \right] =$$
$$= H_n \int f'\left\{V_n(X), Y\right\} dT_n +$$
$$+ \frac{1}{2} H_n \iint f''\left\{V(X), Y, Z\right\} dT_n(Y)\, dT_n(Z).$$

Ici $T_n = S_n - V_n$ et $V(X)$ désigne une distribution sur le « segment » qui va de $S_n(X)$ à $V_n(X)$.

Posons:

$$A_n = H_n \left[ f\{S_n(X)\} - f\{V_n(X)\} \right],$$
$$B_n = H_n \int f'\{V_n(X), Y\}\, dT_n(Y).$$

On a:

$$A_n - B_n = \frac{1}{2} H_n \iint f''\{V(X), Y, Z\}\, dT_n(Y)\, dT_n(Z).$$

Selon le lemme établi, nous montrerons,

1) que la distribution $Q_n(u)$ de $B_n$ tend vers la Gaussienne,

2) que l'espérance mathématique de la valeur absolue de la différence $A_n - B_n$ tend vers zéro si $n$ augmente au-delà de toute valeur.

Nous déduirons, en conséquence, que la distribution $P_n(u)$ de $A_n$ tend aussi vers la Gaussienne, ce qui démontrera le théorème.

**5.3.2.** ***Démonstration de la première proposition***

En tenant compte de la définition de $T_n$, $B_n$ peut s'écrire:

$$B_n = H_n \int f'\{V_n(X), Y\}\, dS_n(Y) - H_n \int f'\{V_n(X), Y\}\, dV_n(Y).$$

Suivant la définition de l'intégrale

$$\int \psi\, dS_n$$

et en remarquant que

$$V_n = \frac{1}{n} \sum_{\nu=1}^{n} V'_\nu$$

on a:

$$B_n = \frac{H_n}{n} \sum_{\nu=1}^{n} f'\{V_n(X), X_\nu\} - \frac{H_n}{n} \sum_{\nu=1}^{n} \int f'\{V_n(X), Y\}\, dV'_\nu(Y).$$

Ici, pour $n$ fixe, $X_1, X_2, ..., X_n$ désignent les résultats immédiats de $n$ épreuves effectuées sur les $n$ premiers collectifs donnés $C'_1, C'_2, ..., C'_n$.

Imaginons que dans le collectif $C'_\nu$ au résultat $X_\mu$, nous attachions la valeur $f'\{V_n(X), X_\nu\}$ et que nous fassions de même pour tous les collectifs. Nous aurons ainsi une suite de collectifs $C'_1, C'_2, ..., C'_n$ ayant pour caractères distinctifs les valeurs respectives

$$f'\{V_n(X), X_1\}, f'\{V_n(X), X_2\}, ..., f'\{V_n(X), X_n\}$$

En ce cas, dans l'expression de $B_n$, la première somme est une somme de $n$ variables aléatoires, chacune à une dimension, et chaque terme de la seconde somme est égal à la valeur moyenne de ces variables.

Posons pour abréger l'écriture:

$$f_\nu = f\{V_n(X), X_\nu\}$$

et en même temps servons nous de la notation $a_\nu$ introduite dans la condition 2 du théorème; nous aurons

$$B_n = \frac{H_n}{n} \sum_{\nu=1}^{n} (f_\nu - a_\nu).$$

A part le facteur $\frac{H_n}{n}$ nous voyons ici, une somme de $n$ variables aléatoires dont les valeurs moyennes s'annulent. La seconde hypothèse sur les fonctions $f\{V(X)\}$ garantit que la distribution de cette somme tend vers la distribution de Gauss. D'après les notations, la dispersion du v$^{\text{ème}}$ de nos variables est égale à $r_\nu^2$, donc la dispersion de leur somme égale

$$s_n^2 = r_1^2 + r_2^2 + ... + r_n^2;$$

enfin la dispersion de $B_n$ est égale à $\left(\frac{H_n}{n}\right)^2 . s_n^2 = \frac{1}{2}$ en vertu de la définition de $H_n$. Donc:

si $Q_n(u)$ désigne la probabilité pour que $B_n$ ne surpasse pas $u$ , on a

$$\lim_{n\to\infty} Q_n(u) = \phi(u) = \frac{1}{\sqrt{\pi}} \int^{u} e^{-t^2} dt$$

**5.3.3.** ***Démonstration de la deuxième proposition***

On a:

$$A_n - B_n = \frac{1}{2} H_n \iint f''\{V(X), Y, Z\} \, dT_n(Y) \, dT_n(Z).$$

En vertu de la troisième hypothèse sur les fonctions $f\{V(X)\}$ on a

$$|A_n - B_n| \le \frac{H_n}{2} \int \psi(X) \left[S_n(X) - V_n(X)\right]^2 dX.$$

Il faut trouver l'espérance mathématique du second membre. Désignons l'intégrale par $J$ , on a:

$$|A_n - B_n| \le \frac{H_n}{2} J$$

et suivant une formule mentionnée au § 2:

$$E_n\{|A_n - B_n|\} \le \frac{H_n}{2n} E_n\{J\} \le \frac{H_n}{2n} \int \psi(X) V_n(X) \left[1 - V_n(X)\right] dX.$$

Mais $\dfrac{H_n}{2n} = \dfrac{1}{2s_n\sqrt{2}}$, d'où

$$E_n\{|A_n - B_n|\} \le \frac{1}{2s_n\sqrt{2}} \int \psi(X) \left[1 - V_n(X)\right] V_n(X) \, dX.$$

En passant à la limite et en vertu de l'hypothèse 3, on a

$$\lim_{n\to\infty} E_n\{|A_n - B_n|\} = 0.$$

Les deux conditions du lemme sont remplies; par suite, si $P_n(u)$ désigne la probabilité pour que la valeur de $A_n$ définie par

$$H_n\left[ f\{ S_n(X)\} - f\{V_n(X)\}\right]$$

ne surpasse pas $u$, on a que :

$$\lim_{n\to\infty} P_n(u) = \phi(u) = \frac{1}{\sqrt{\pi}} \int^u e^{-t^2} dt.$$

Le théorème se trouve donc démontré.

**5.4.** ***Remarques sur les fonctions $f\{V(X)\}$ satisfaisant à nos conditions.***

1. Pour les fonctions linéaires $\int \alpha(X)\, dV(X)$ la question ne se pose pas: la deuxième dérivée s'annule, la première dérivée est égale à $\alpha(X)$ et est indépendante de $V$. Les seules conditions à remplir, sont les conditions nécessaires pour que le théorème classique subsiste pour les variables aléatoires soumises aux lois de distribution $V_\nu'(X)$.

2. Pour les fonctions non linéaires, le type qui se présente le plus souvent est celui de la forme

$$f = F(A,B,C,\ldots)$$

avec

$$A = \int \alpha(X)\, dV(X),\ \ B = \int \beta(X)\, d(V - V_1)(X),\ldots$$

En prenant pour $V(X)$ des distributions restreintes aux $V(X)$ à partir d'un certain $n$, nous pourrons supposer que les dérivées du premier et du second ordre de $F$ par rapport à $A,B,C,\ldots$ soient bornées. Mais si l'on se rappelle que la première dérivée de $f\{V(X)\}$ définie ci-dessus, se compose de termes de la forme

$$\frac{\partial F}{\partial A}\alpha(y),\ \frac{\partial F}{\partial B}\beta(y),\ldots$$

la deuxième condition est remplie si chacune des fonctions $\alpha,\beta$, … satisfait aux conditions du théorème classique au sens indiqué ci-dessus.

Pour la troisième hypothèse, remarquons que la deuxième dérivée dépend d'expressions de la forme

$$\frac{\partial^2 F}{\partial A^2}\alpha(Y)\alpha(Z),\ \frac{\partial^2 F}{\partial A\partial B}\alpha(Y)\beta(Z),\dots$$

Comme les $\frac{\partial^2 F}{\partial A^2}$ restent bornés, nous aurons des expressions de la forme

$$\iint \alpha(Y)\beta(Z)\,dT_n(Y)\,dT_n(Z)=\int \alpha(Y)\,dT_n(Y)\cdot\int \beta(Z)\,dT_n(Z).$$

Utilisons maintenant la formule d'intégration par parties suivante

$$\begin{aligned}\iint \alpha(x,y)\,dT_n(x,y)= & \\ &=\iint \alpha_{xy}T\,dx\,dy-\\ &-\int \alpha_x(x,\infty)\,T(x,\infty)\,dx-\\ &-\int \alpha_y(\infty,y)\,T(\infty,y)\,dy\end{aligned}$$

que l'on peut établir facilement. Nous supposons d'abord les collectifs à deux dimensions.

Les distributions $V_\nu'$ faisant partie de l'ensemble convexe $J$ sont choisies de façon que les intégrales ci-dessus soient convergentes.

Dès lors en vertu de la convergence des intégrales on peut écrire

$$\int |\alpha_x(x,\infty)|\,|T(x,\infty)|\,dx\le\iint |\alpha_x(x,y)|\,|T(x,y)|\,dx\,dy.$$

D'une façon identique on a

$$\int |\alpha_y(\infty,y)|\,|T(\infty,y)|\,dy\le\iint |\alpha_y(x,y)|\,|T(x,y)|\,dx\,dy.$$

Par suite en égard de notre formule d'intégration par parties on a:

$$\left|\iint \alpha(x,y)\, dT_n(x,y)\right| \leq$$

$$\leq \iint \left( \left|\alpha_{xy}(x,y)\right| + \left|\alpha_x(x,y)\right| + \left|\alpha_y(x,y)\right| \right) \left|T(x,y)\right| dx\,dy.$$

Supposons maintenant que les dérivées $\alpha_{xy}(x,y), \alpha_x(x,y), \alpha_y(x,y)$ de $\alpha, \beta_{xy}, \beta_x, \beta_y$ de $\beta$ et ainsi de suite aient une même majorante $\psi_1(x,y)$ en valeur absolue.

Appliquons l'inégalité de Schwarz pour une fonction positive $\psi(x,y)$; nous aurons:

$$\left|\iint \alpha(x,y)\, dT_n(x,y)\right|^2 \leq 9 \iint \frac{\psi_1^2(x,y)}{\psi(x,y)} dx\,dy \cdot \iint \psi(x,y)\, T_n^2(x,y)\, dx\,dy$$

que l'on peut écrire aussi en tenant compte de ce que

$$\iint \alpha(Y)\beta(Z)\, dT_n(Y)\, dT_n(Z) = \int \alpha(Y)\, dT_n(Y) \cdot \int \beta(Z)\, dT_n(Z):$$

$$\left|\iint \alpha(Y)\beta(Z)\, dT_n(Y)\, dT_n(Z)\right| \leq$$

$$\leq 9 \iint \frac{\psi_1^2(x,y)}{\psi(x,y)} dx\,dy \cdot \iint \psi(x,y)\, T_n^2(x,y)\, dx\,dy.$$

Par suite pour que la troisième hypothèse soit remplie, il suffit que l'intégrale

$$\iint \frac{\psi_1^2(x,y)}{\psi(x,y)} dx\,dy$$

soit convergente pour toute fonction $\psi(x,y)$ pour laquelle

$$\lim_{n\to\infty} \frac{1}{s_n} \iint \psi(x,y) V_n(x,y) \left[1 - V_n(x,y)\right] dx\,dy = 0$$

$s_n$ étant défini dans l'énoncé du théorème.

Un raisonnement analogue pour le cas de plus de deux dimensions, nous permet de montrer que

$$\left|\iint\dots\int \alpha(x_1,x_2,...,x_k)\,dT(x_1,x_2,...,x_k)\right| \leq$$

$$\leq \iint\dots\int\Big[\left|\alpha_{x_1,x_2,...,x_k}\right| + \dots$$

$$+\left|\alpha_{x_1,x_2,...,x_{k-1}}\right| + \dots \left|\alpha_{x_k}\right|\Big]\left|T(x_1,x_2,...,x_k)\right| dx_1 dx_2 ... dx_k .$$

l’expression entre crochets, désigne la somme de toute les dérivées de $\alpha$ par rapport à chaque variable jusqu'à l'ordre $k$ inclus prises en valeur absolue. En supposant que les dérivées des fonctions $\alpha,\beta,\gamma,...$ jusqu'à l'ordre $k$ inclus aient une même majorante $\psi_1(X)=\psi_1(x_1,x_2,...,x_k)$ en valeur absolue, l'inégalité de Schwarz appliquée pour une fonction positive $\psi(x_1,x_2,...,x_k)$ conduit au résultat que

$$\left|\iint \alpha(Y)\beta(Z)\,dT_n(Y)\,dT_n(Z)\right| \leq$$

$$\leq K\cdot\iint\dots\int \frac{\psi_1^2(x_1,...,x_k)}{\psi(x_1,...,x_k)}\,dx_1...dx_k \cdot \iint\dots\int \psi T_n^2 dx_1...dx_k$$

$K$ étant constante.

Dès lors, si l’intégrale

$$\iint\dots\int \frac{\psi_1^2}{\psi}\,dX$$

converge pour toute fonction $\psi(X)$ telle que

$$\lim_{n\to\infty}\frac{1}{s_n}\int\dots\int \psi(X)V_n(X)\big[1-V_n(X)\big]dX = 0$$

la troisième hypothèse est remplie.

**5.5.** ***Cas du coefficient de corrélation et de fonctions quelconques des moments.***

**5.5.1.** Supposons les collectifs à deux dimensions. Le coefficient de corrélation $\Gamma$ se définit par

$$\Gamma = \frac{M_{11}}{\sqrt{M_{20}M_{02}}} = \frac{\iint (x-a)(y-b)\,dV(x,y)}{\iint (x-a)^2\,dV(x,y)\cdot\iint (y-b)^2\,dV(x,y)}$$

avec

$$a = \iint x dV(x,y) \qquad b = \iint y dV(x,y)$$

Afin de calculer les dérivées statistiques de ce coefficient, écrivons-le ainsi:

$$\Gamma = \frac{\iint xy dV - ab}{\sqrt{\left[\iint x^2 dV - a^2\right]\left[\iint y^2 dV - b^2\right]}}$$

Posons

$$A = \iint xy dV,\ a = B = \iint x dV,\ b = C = \iint y dV,\ D = \iint x^2 dV,\ E = \iint y^2 dV$$

et l'on a la fonction:

$$\Gamma = \frac{A - BC}{\sqrt{\left(D - B^2\right)\left(E - C^2\right)}}$$

où $A, B, \ldots$ sont des fonctions statistiques linéaires. Pour avoir les dérivées de $\Gamma$, nous appliquerons les règles établies pour les fonctions du genre $F(A, B, C, \ldots)$. Nous désignerons par $y_1, y_2$ les composantes du point $Y$ et par $z_1, z_2$ celles du point $Z$ qui figurent dans la première et dans la seconde dérivée. Pour déterminer les dérivées statistiques, il faut calculer les dérivées de $\Gamma$ par rapport à $A, B, C....$

On obtient ainsi par des calculs faciles:

$$f'\{V(X),Y\} =$$
$$= \frac{\Gamma}{A} y_1 y_2 + \frac{\Gamma}{A'D'}\left(A'B - D'C\right) y_1 +$$
$$+ \frac{\Gamma}{A'E'}\left(A'C - E'B\right) y_2 - \frac{1}{2}\frac{\Gamma}{D'} y_1^2 - \frac{1}{2}\frac{\Gamma}{E'} y_2^2 .$$

Ici

$$A' = A - BC = M_{11},\ D' = D - B^2 = M_{20},\ E' = E - C^2 = M_{02} .$$

Avec la notation des moments, en tenant compte que

$$\Gamma = \frac{M_{11}}{\sqrt{M_{20}M_{02}}}$$ , on obtient de définitive:

$$f'\{V(X),Y\} =$$
$$= \frac{M_{11}}{\sqrt{M_{02}M_{20}}}\left[\frac{1}{M_{11}} y_1 y_2 + \left(\frac{a}{M_{20}} - \frac{b}{M_{11}}\right) y_1 + \right.$$
$$\left. + \left(\frac{b}{M_{02}} - \frac{a}{M_{11}}\right) y_2 - \frac{1}{2}\frac{1}{M_{20}} y_1^2 - \frac{1}{2}\frac{1}{M_{02}} y_2^2 \right].$$

Pour la dérivée seconde on obtient par des calculs :

$$f''\left\{V(X),Y,Z\right\}=$$

$$=\frac{M_{11}}{\sqrt{M_{02}M_{20}}}\left[\frac{2a}{M_{11}M_{20}}y_1y_2z_1+\frac{2b}{M_{11}M_{02}}y_1y_2z_2-\frac{1}{M_{11}M_{20}}y_1y_2z_1^2+\right.$$

$$+2\left(-\frac{1}{M_{11}}-\frac{a^2}{M_{11}M_{20}}-\frac{b^2}{M_{11}M_{02}}+\frac{ab}{M_{02}M_{20}}\right)y_2z_1+$$

$$+\frac{1}{2}\frac{1}{M_{02}M_{20}}y_2^2z_1^2-\frac{1}{M_{11}M_{02}}y_1y_2z_2^2+$$

$$+\left(\frac{1}{M_{20}}-\frac{2ab}{M_{11}M_{20}}+\frac{3a^2}{M_{20}^2}\right)y_1z_1+$$

$$+2\left(-\frac{1}{M_{11}}-\frac{b^2}{M_{11}M_{02}}-\frac{a^2}{M_{11}M_{20}}+\frac{ab}{M_{20}M_{02}}\right)\left(y_1z_2+y_2z_1\right)$$

$$+\left(\frac{b}{M_{11}M_{20}}-\frac{3a}{M_{20}^2}\right)y_1z_1^2+\left(\frac{b}{M_{11}M_{02}}-\frac{a}{M_{20}M_{02}}\right)y_1z_2^2+$$

$$+\left(\frac{-2ab}{M_{11}M_{02}}+\frac{1}{M_{02}}+\frac{3b^2}{M_{02}^2}\right)y_2z_2+$$

$$+\left(\frac{a}{M_{11}M_{20}}-\frac{b}{M_{20}M_{02}}\right)y_2z_1^2+\left(\frac{a}{M_{11}M_{02}}-\frac{3b}{M_{02}^2}\right)y_2z_2^2+$$

$$\left.+\frac{3}{4}\frac{1}{M_{20}^2}y_1^2z_1^2+\frac{1}{2}\frac{1}{M_{02}M_{20}}\left(y_1^2z_2^2+y_2^2z_1^2\right)+\frac{3}{4}\frac{1}{M_{02}^2}y_2^2z_2^2\right].$$

Cherchons à appliquer les conditions de notre théorème. Remarquons que suivant les notations, les fonctions $\alpha,\beta,\gamma,\ldots$ paraissant dans les intégrales $A,B,C,\ldots$ qui figurent dans $f\left\{V(X)\right\}=F(A,B,C,\ldots)$ sont:

$\alpha(x,y)=x,\ \beta(x,y)=y,\ \gamma(x,y)=xy$ (ou aussi bien $(x-a)(y-b)$)

$\delta(x,y) = x^2$ (ou $(x-a)^2$), $\varepsilon(x,y) = y^2$ (ou $(y-b)^2$) (*a* et *b* sont des valeurs restreintes à un intervalle fini).

Les dérivées partielles sont les suivantes:

$$\alpha_x = 1,\ \alpha_y = 0,\ \alpha_{xy} = 0\ ; \beta_x = 0\ ,\ \beta_y = 1\ ,\ \beta_{xy} = 0\ ;$$

$$\gamma_x = y\ ,\ \gamma_y = x\ ,\ \gamma_{xy} = 1\ ;\ \delta_x = 2x\ ,\ \delta_y = 0\ ,\ \delta_{xy} = 0\ ,\ldots$$

Il est facile de voir que toutes ces dérivées ont la majorante en valeur absolue:

$$\psi_1(x,y) = c + 2\sqrt{x^2+y^2} \text{ avec } c > 1.$$

L'intégrale $\iint \frac{\psi_1^2(x,y)}{\psi(x,y)}\,dxdy$

dont il est question sera donc convergente si $\psi(x,y)$ augmente vers l'infini comme $\left(\sqrt{x^2+y^2}\right)^{4+\varepsilon}$ $(\varepsilon > 0)$ quand *x* et *y* surpassent simultanément toute quantité donnée.

Dès lors, si les distributions données $V'(x,y)$ remplacent la condition que le produit

$$\left|\sqrt{x^2+y^2}\right|^{6+\varepsilon} V_\nu'(x,y)$$

pour $x < -N$, $y < -N'$ ($N$ et $N' > 0$ et très grands) et le produit

$$\left|\sqrt{x^2+y^2}\right|^{6+\varepsilon} \left[1 - V_\nu'(x,y)\right]$$

pour $x > N$ et $y > N'$ restent inférieurs à un nombre indépendant de $\nu$, les intégrales

$$\iint \psi(x,y) V_n(x,y)\left[1 - V_n(x,y)\right] dxdy$$

seront convergentes uniformément par rapport à *n* et ainsi la troisième hypothèse est vérifiée. Ainsi dans ces conditions on aura le résultat intéressant que *la distribution du coefficient de corrélation tendra vers la Gaussienne.*

**5.5.2.** Dans les cas de distributions de plus de deux dimensions et si, en général, la fonction $f$ considérée dépend de plusieurs moments, $M_{\nu_1 \ldots \nu_k}$ jusqu'à l'ordre $m$ (le plus grand parmi les $\nu$ ) les dérivées des fonctions $\alpha, \beta$ … de $A, B$ … auront la majorante, en valeur absolue

$$c + c'|X|^{m-1}$$

ou $c$ et $c'$ sont deux constantes $(c > 1)$ et $|X|$

la longueur du vecteur $X : |X| = \sqrt{x_1^2 + x_2^2 + \ldots + x_k^2}.$

Dès lors pour que l'intégrale $\int \frac{\psi_1^2}{\psi} dX$ soit convergente, il suffit que $\psi(X)$ augmente comme $|X|^{2m-2+k+\varepsilon}$

pour $x_i > N_i \qquad (i = 1, 2, \ldots, k)$

$N_i$ étant un nombre positif très grand.

Mais en considérant l'intégrale $\int \psi(X) V_n(X)\left[1 - V_n(X)\right] dX$ on voit qu'elle sera convergente, si les distributions données $V'_\nu(X)$ remplissent les conditions que

$V'_\nu(X)$ pour $x_i < -N_i$ $\qquad$ $N_i$ est un nombre positif très grand.

$1 - V'_\nu(X)$ pour $x_i > N_i$

$$(i = 1, 2, \ldots, k)$$

restent inférieurs à $\lambda |X|^{-2m+2-2k-\varepsilon}$, $\lambda$ étant une constante indépendante de $\nu$ $(\varepsilon > 0)$.

Si l'on rappelle que la seconde hypothèse sur les fonctions $f\{V(X)\}$ demande que des intégrales de la forme $\int \alpha^{2m+\varepsilon} dV'_\nu$ soient convergentes, on constate que les fonctions $V'_\nu$ doivent remplir les conditions que

$$V'_\nu(X) \text{ pour } x_i < -N_i$$

$$1 - V'_\nu(X) \text{ pour } x_i > N_i$$

$$N_i >> 0$$
$$(i = 1, 2, ..., k)$$

restent inférieurs, indépendamment de $v$, à

$$\lambda |X|^{-2m-k-\varepsilon} \qquad (\varepsilon > 0)$$

$\lambda$ étant une constante.

***Remarque:***

Dans la considération de la convergence des intégrales multiples, on s'est basé sur le fait que l'intégrale généralisée

$$\iint \dots \int \frac{dx_1 dx_2 \dots dx_k}{\left(x_1^2 + x_2^2 + \dots + x_k^2\right)^m}$$

converge si $m > \frac{1}{2}k$.

*****